\documentclass[10pt, english]{amsart}

\usepackage{amsmath,amssymb,enumerate}

\usepackage[T1]{fontenc}
\usepackage[all,cmtip]{xy}

\usepackage{babel}
\usepackage{amstext}
\usepackage{amsmath}
\usepackage{amsfonts}
\usepackage{latexsym}
\usepackage{ifthen}

\usepackage{xypic}
\xyoption{all}
\newtheorem{lemma1}{}[section]

\newenvironment{theorem}{\begin{lemma1}{\bf Theorem.}}{\end{lemma1}}

\newenvironment{proposition}{\begin{lemma1}{\bf Proposition.}}{\end{lemma1}}
\newenvironment{corollary}{\begin{lemma1}{\bf Corollary.}}{\end{lemma1}}

\newenvironment{definition}{\begin{lemma1}{\bf Definition.}}{\end{lemma1}}

\newenvironment{problem}{\begin{lemma1}{\bf Problem.}}{\end{lemma1}}

\newenvironment{the local obstruction - setup}{\begin{lemma1}{\bf The local obstruction - setup.}}{\end{lemma1}}

\newenvironment{remark*}{{\bf Remark.}}{}
\newenvironment{remarks*}{{\bf Remarks.}}{}
\newenvironment{example*}{{\bf Example.}}{}
\newenvironment{assumption*}{{\bf Assumption.}}{}

\newcommand{\PP}{\ensuremath{\mathbb{P}}}

\usepackage{eurosym}

\makeatletter
\ifnum\@ptsize=0  \fi \ifnum\@ptsize=2  \fi 

\newcommand\sO{{\mathcal O}}

\usepackage{xcolor}
\usepackage{hyperref}

\author{Thomas Peternell}

\address{Thomas Peternell, Mathematisches Institut, Universit\"at Bayreuth, 95440 Bayreuth, 
Germany}
\email{thomas.peternell@uni-bayreuth.de}

\title{Compactifications of $\mathbb C^n$  and the complex projective space } 

\begin{document}

\maketitle

\begin{abstract}We prove that the complex projective space is the only compact K\"ahler manifold compactiying $\mathbb C^n$ by a smooth connected hypersurface, provided $n$ is even. In the odd-dimensional 
case, partial results are given.  The case $n \equiv 1$ mod $4$ has now been settled by Ping Li. 
\end{abstract}


\section{Introduction} 

In his famous problem list \cite{Hi54}, Hirzebruch asked to classifiy all (smooth) compactifications $X$ of $\mathbb C^n$ with $b_2(X) =1$. The condition on $b_2(X)$ is equivalent to saying that the divisor at $\infty$ is
irreducible. 
For an overview on this problem,
see \cite{PS89}. In any case, without additional assumptions, such a classification seems to be possible only in low dimensions. 

A particularly interesting case is when $X$ is K\"ahler and the divisor at infinity
is smooth. A folklore conjecture states that then $X$ must be projective space. In this paper, we confirm this conjecture if $X$ is of even dimension. More generally, we show

\begin{theorem} \label{thmA} Let $X$ be a compact K\"ahler manifold of dimension $n$ and $Y \subset X$ a smooth connected hypersurface such that  
$$ H_p(Y,\mathbb Z) \to H_p(X,\mathbb Z)$$
is bijective for all $0 \leq p \leq 2n-2$. 
\begin{enumerate} 
\item If $n$ is even, then  $X \simeq \mathbb P_n$ and $Y$ is a hyperplane. 
\item If $n$ is odd, then either  $X \simeq \mathbb P_n$ and $Y$ is a hyperplane or $X$ is a Fano manifold with Picard number one, of index $\frac{n+1}{2}$,
with $Y \in \vert \sO_X(1) \vert $. Further, $c_{n-1}(X) = n \cdot n+1$. 
\item Let $\sO_X(1) $ be the ample generator of ${\rm Pic}(X) \simeq \mathbb Z$. If $n$ is odd and if $h^0(X,\sO_X(1) \geq 2$, then $X \simeq \mathbb P_n$ and $Y$ is a hyperplane.
\end{enumerate} 
 \end{theorem}
 
 The case $n \equiv 1$ mod $4$ has now been settled by Ping Li \cite{Li25}. 
 
 The condition in Item (3) can be rephrased by $\chi(X,\sO_X(1)) \ne 1$ or $\chi(Y,\sO_Y(1)) \ne 0$. If $n = 7$ or $n = 9$, condition (c) can be verified. In Section 3 we present some approaches to the 
 odd-dimensonal case. In particular, we set up a system of equations in terms of the Chern classes $c_j(X)$ and $c_j(Y)$.  If this set has no  integer solution for $c_1(X) = \frac{n+1}{2}$, then Theorem \ref{thmA} 
 is established. 
 
 \begin{corollary} Let $X$  be a compact K\"ahler manifold of even dimension $n$ and $Y \subset X$ a smooth connected hypersurface such that  
 $X \setminus  Y$ is biholomorphic to $\mathbb C^n$. Then $X \simeq \mathbb P_n$ and $Y$ is a hyperplane. 
 \end{corollary} 
 
 This result was shown independently without assumption on the dimension by completely different methods by Chi Li and Zhengyi Zhou, \cite{LZ24}.

Theorem \ref{thmA}(a) has been proved  for  any $n \leq 5$ by van de Ven \cite{vdV62} and for $n \leq 6$  by Fujita \cite{Fu80}.

For important discussion and remarks I would like to thank Daniel Barlet, Fr\'ed\'eric Campana, Baohua Fu, Andreas H\"oring and Mihai P\v{a}un.

\section{Proof of Theorem \ref{thmA} }

We fix $X$ and $Y$ as in Theorem \ref{thmA} and  first collect some basic properties of $X$ and $Y$ and refer to  Fujita \cite{Fu80} and Sommese \cite{So76}.

\begin{proposition}  \label{BP}
\begin{enumerate}
\item The cohomology ring $H^*(X,\mathbb C)$ is isomorphic as graded ring to $H^*(\mathbb P_n,\mathbb C)$.
In particular  $h^{p,q}(X) = h^{p,q}(\PP_n) $ for all $p,q$.
\item  The cohomology ring $H^*(Y,\mathbb C)$ is isomorphic as graded ring to $H^*(\mathbb P_{n-1},\mathbb C)$.
In particular  $h^{p,q}(Y) = h^{p,q}(\PP_{n-1}) $ for all $p,q$.
\item The restrictions $H^q(X,\mathbb C) \to H^q(Y,\mathbb C)$ are bijective for $0 \leq q \leq 2n-2$.
\item $X$ is a Fano manifold. Further ${\rm Pic} (X) \simeq \mathbb Z$, with ample generator $\mathcal O_X(Y)$.
\end{enumerate}
\end{proposition} 

Thus we may regard all Chern numbers of $X$ and $Y$ as numbers and intersection is just multiplication.  Further let $r$ denote the index of $X$, so that $-K_X = \sO_X(rY)$. In other words,
$$ c_1(X) = r.$$
Note that it is not necessary to assume $Y$ to be ample, since ${\rm Pic}(X) \simeq \mathbb Z$ holds actually for any smooth compactification with $b_2 = 1$ as well as (1) and (2) in the proposition. 


By Proposition \ref{BP} and Corollary 2.5 of Libgober-Wood \cite{LW90}, we have

\begin{proposition} \label{prop1} 
\begin{equation} \label{eq1}  r  \cdot  c_{n-1}(X)  = c_1(X) \cdot c_{n-1}(X) = \frac{1}{2} n (n+1)^2 \end{equation} 
 and 
\begin{equation} \label{eq2}  (r-1)  \cdot c_{n-2}(Y) =  c_1(Y) \cdot c_{n-1}(Y) = \frac{1}{2} (n-1) n^2.\end{equation} 
 \end{proposition}

 We now start the proof of Theorem \ref{thmA} and observe first that $r \ne 1$ by Equation (\ref{eq2}).
 The tangent bundle sequence
 $$ 0 \to T_Y \to T_X \vert Y \to N_{Y/X} \to 0$$ yields  (in terms of numbers!) 
 $$ c_{n-1}(X) =  c_{n-1}(T_X \vert Y) = c_{n-1}(Y) + c_{n-2}(Y) \cdot c_1(N_{Y/X}) =  c_{n-1}(Y) + c_{n-2}(Y).$$
 Since $$ c_{n-1}(Y) = \chi_{\rm top} (Y) = \chi_{\rm top}(\mathbb P_{n-1}) = n,$$ it follows 
 
 \begin{equation} \label{eq3} c_{n-1}(X) = n+ c_{n-2}(Y).\end{equation} 
 Replacing $c_{n-1}(X)$ by $n+ c_{n-2}(Y)$ in Equation (\ref{eq1}) and putting in 
 $$ c_{n-2}(Y) = \frac{1}{2(r-1)} (n-1)n^2$$
 by virtue of Equation (\ref{eq2}), we obtain

 \begin{equation} \label{eq4}  r\big(1 + \frac{1}{2(r-1)} (n-1) n\big) = \frac{1}{2} (n+1)^2. \end{equation} 
 
Fixing $r$, we obtain a quadratic equation for $n$ with solutions
 $$ n = r-1 $$ and
 $$ n = 2r-1.$$ 
 In the first case $X \simeq \mathbb P_n$ proving (a)
 
In case $ n = 2r-1$, $n$ is odd. If $X$ is not projective space, then 
 $$ c_1(X) = \frac{n+1}{2}  = \frac{1}{2} c_1(\mathbb P_n)$$
 and
 $$ c_{n-1}(X) = n(n+1) = 2 c_{n-1}(\mathbb P_n),$$ establishing (b). 
 
 So assume finally that $h^0(X,\sO_X)) \geq 2$ and pick $Z \in \vert  \sO_X(1) \vert $ general. Then $Z $ is smooth, different from $Y$. Since 
 $Y \cdot Z$ has class one in $H^4(X,\mathbb Q)$, it follows from \cite[Prop.7.2]{Fu84} that $Y \cap Z$ is smooth. Further, by Mayer-Vietoris, the inclusions
 $$H_q(Y\cap Z,\mathbb C) \to H_q(Y,\mathbb C)$$ 
 are bijective for $0 \leq q \leq 2n-4$.
 Since $\dim Y$ is even, it follows from Part (a) that $Y \simeq \mathbb P_{n-1}$. Hence $X \simeq \mathbb P_n$ using again Equations (\ref{eq1}) and (\ref{eq2}). 
 
 \section{Approaches to the odd-dimensional case} 
 
 We extract the non-vanishing in Theorem \ref{thmA}(c) as follows. Recalling that the restriction map
 $$ H^0(X,\sO_X(1)) \to H^0(Y,\sO_Y(1)) $$
 is surjective due to the Kodaira vanishing theorem, we are reduced to the following 
 
 \begin{problem} \label{Pr2} Let $Y$ be a Fano manifold of even dimension $d$ and index $r = \frac{d}{2}$. Assume $\rho(Y) = 1$. 
 Let $\sO_Y(1)$ be ample such that $-K_Y = \sO_Y(r)$ and assume that $c_1(\sO_Y(1))^n  = 1$. Is
 $$ \dim  H^0(Y,\sO_Y(1)) \ne 0?$$ 
 \end{problem} 
 
 Of course, we may assume much more in our setting: $Y$ has the same Hodge numbers as $\mathbb P_d$, hence $c_{d-1}(Y) = d(d+1)$ and the cohomology ring is that of $\mathbb P_d$.

 One checks by hand that this problem has a positive solution for $d \leq 8$. E.g., assume $d = 8$. Consider the Hilbert polynomial
 $$ p(t) = \chi(\sO_Y(t)). $$
 Assuming that $H^0(Y,\sO_Y(1)) = 0$, the polynomial $p$ has zeroes at $t = 1,-1,-2,-3, -5.$ Further $p(0) = 1 = p(-4)$. Taking into account that the coefficient 
 of $t^7$ is $\frac{4}{7!}$ by Riemann-Roch, and that $p(t) = p(-4-t)$ for all $t$ (by Serre duality), one derives a contradiction by explicit calculation.

 This yields Theorem \ref{thmA}(a) in dimenson $7$ and $9$.

 \begin{definition} 
Let $Z$ be any compact manifold of dimension $n$. Set 
 $$ \chi_Z(t) := \sum_{p=0}^n \chi_p(Z) t^p $$
 where 
 $$ \chi_p(Z) = \sum_{q=0}^n (-1^q) h^q(\Omega_Z^p). $$
 Set further 
 $$ A_k(Z) = \frac{1}{(2k)!} \chi_Z^{(2k)}(-1). $$ 
 \end{definition}

 Note that $A_0(Z) = c_n(Z) $ and that 
 $A_1(Z) $ is a linear combination of $c_n(Z)$ and $c_1(Z) \cdot c_{n-1}(Z)$. 
 Further, $A_2(Z)$ is a linear combination of $A_0(Z) $ and $A_1(Z)$ and of 
 $c_{n-2}(Z) $ and $c_{n-3}(Z)$, see \cite[Prop.2.3]{LW90} and \cite[p.145]{Sa96}.
 Basically, the new term in $A_2(Z)$ is 
 $$ \big(c_1^2+3c_2\big) \cdot c_{n-2} - \big(c_1^3-3c_1c_2+3c_3\big) \cdot c_{n-3}.$$
 
 The point is of course

 \begin{proposition} Let $X$ and $Y$ be compact complex manifolds of dimension $n \geq 3$. If $X$ and $Y$ have the same 
 Hodge numbers, then 
 $$ A_k(X) = A_k(Y) $$ for all $k$. 
 \end{proposition} 
 
 Hence:

 \begin{proposition}  Under the assumptions of Theorem \ref{thmA}, the following holds for all $0 \leq k \leq n$ (resp. $0 \leq k \leq n-1)$
 \begin{equation} A_k(X) = A_k(\mathbb P_n)  \end{equation}
 and
  \begin{equation} A_k(Y) = A_k(\mathbb P_{n-1})   \end{equation} 
  \end{proposition} 
  
  \begin{proof} Apply Proposition \ref{BP}. 
  \end{proof} 
  
 Thus we obtain a system of equations 
 \begin{equation} \label{EQ1}  A_k(X) = A_k(\mathbb P_n), \  0 \leq k \leq n, \end{equation} 
  \begin{equation} \label{EQ2}  A_k(Y) = A_k(\mathbb P_{n-1}),  \ 0 \leq k \leq n-1, \end{equation} 
  and 
  \begin{equation} \label{EQ3} c_j(X) = c_j(Y) + c_{j-1}(Y),  \ 1 \leq j \leq n. \end{equation} 
  
  Actually, (\ref{EQ1}) gives $\frac{n+1}{2}$ independent equations and so gives (\ref{EQ2}), see \cite{LW90}, p.142. 
  The hope (strong version) is now that this system of equation has only one integer solution, namely $c_j(X) = c_j(\mathbb P_n) = \binom{n+1}{j}$ for all $j$ 
  (and thus $c_j(Y) = c_j(\mathbb P_{n-1})$). It would be however be sufficient to prove a weak version, namely that there is no integer solution with $c_1(X) = \frac{n+1}{2}$. 
  
Recall that the equations (\ref{EQ1}) and (\ref{EQ2}) for $k = 0$ and $k = 1$ simply give 
$$ c_n(X) = c_n(\mathbb P_n) = n+1;$$
$$ c_{n-1}(Y) = c_{n-1}(\mathbb P_{n-1}) = n; $$
$$ c_1(X) \cdot c_{n-1}(X) = c_1(\mathbb P_n) \cdot  c_{n-1}(\mathbb P_n) = \frac{1}{2} n (n+1)^2;$$
$$ c_1(Y) \cdot c_{n-2}(Y) = c_1(\mathbb P_{n-1}) \cdot  c_{n-2}(\mathbb P_{n-2}) = \frac{1}{2} (n-1) n^2.$$
For $n = 5$, we get two more equations 
 $$ \big(c_1^2(X) + 3c_2(X)\big) \cdot c_{3}(X) - \big(c_1^3(X) - 3c_1(X) \cdot c_2(X) + 3c_3(X)\big) \cdot c_{2}(X) =  $$
  $$ =  \big(c_1^2(\mathbb P_5) + 3c_2(\mathbb P_5)\big) \cdot c_{3}(\mathbb P_5) - \big(c_1^3(\mathbb P_5) - 3c_1(\mathbb P_5) \cdot c_2(\mathbb P_5) + 3c_3(\mathbb P_5)\big) \cdot c_{2}(\mathbb P_5)
   $$ and
  $$ \big(c_1^2(Y) + 3c_2(Y)\big) \cdot c_{2}(Y) - \big(c_1^3(Y) - 3c_1(Y) \cdot c_2(Y) + 3c_3(Y)\big) \cdot c_{1}(Y) =  $$
  $$ =  \big(c_1^2(\mathbb P_4) + 3c_2(\mathbb P_4)\big) \cdot c_{2}(\mathbb P_4) - \big(c_1^3(\mathbb P_4) - 3c_1(\mathbb P_4) \cdot c_2(\mathbb P_4) + 3c_3(\mathbb P_4)\big) \cdot c_{2}(\mathbb P_4).$$
  These equations in combination with the equations (\ref{EQ3}) do not have an integer solution in case $c_1(X) = 3$.

 Here are some further informations on the Chern classes of $X$.

 \begin{proposition} $\sum_{k=0}^{n} (-1)^k c_k(X) = (-1)^n.$
 \end{proposition} 
 
 \begin{proof}  By 
 the log version of Hopf's theorem, see \cite[2.1]{CMZ22},
 we have
 $ c_n(\Omega^1_X(\log Y)) = {(-1)}^n.$
 Since $c_k(\sO_Y) = 1$ (here we identify $\sO_Y $ and $j_*(\sO_Y)$, where
 $j: Y \to X$ is the inclusion), the formula follows. 
 
 \end{proof}

By \cite{Li19}, 
  \begin{equation} \label{Li}  c_2(X) \geq \frac{1}{8} (n^2-1). \end{equation}

 \begin{proposition} Suppose we know that $p_1(X) = p_1(\mathbb P_n).$
 Then $X = \mathbb P_n$. 
\end{proposition} 

\begin{proof} 
Assume $X \ne \mathbb P_n$. Then
$$ p_1(X) = \frac{1}{4} (n+1)^2 - 2c_2(X).$$
On the other hand, $p_1(\mathbb P_n) = n+1$. 
Hence
$$ c_2(X) = \frac{n^2}{8} - \frac{n}{4} - \frac{3}{8}.$$ 
This contradicts Li's result, Equation (\ref{Li}).
 
 \end{proof}

\providecommand{\bysame}{\leavevmode\hbox to3em{\hrulefill}\thinspace}
\providecommand{\MR}{\relax\ifhmode\unskip\space\fi MR }
\providecommand{\MRhref}[2]{%
  \href{http://www.ams.org/mathscinet-getitem?mr=#1}{#2}
}
\providecommand{\href}[2]{#2}


\end{document}